\documentstyle[11pt,twoside]{article}
\newtheorem{theorem}{\noindent\bf Theorem}[section]
\newtheorem{proposition}[theorem]{\noindent\bf Proposition}

\newtheorem{remark}[theorem]{\noindent\bf Remark}

\newtheorem{definition}[theorem]{\noindent\bf Definition}

\newcommand{\e}{\hfill\blacksquare}
\input{amssym}
\begin{document}
\date{}
\title{{\Large\bf Invariances of the operator properties of frame multipliers
under perturbations of frames and symbol}}

\author{{\normalsize\sc H. Javanshiri}}

\maketitle
\normalsize

\begin{abstract}
Let $\Phi$ and $\Psi$ be frames for $\cal H$ and let
$M_{m,\Phi,\Psi}$ be a frame multiplier with the symbol $m$. In
this paper, we restrict our investigation to show that the
operator properties of $M_{m,\Phi,\Psi}$ are stable under the
perturbations of $\Phi$, $\Psi$ and $m$. Also, special attention
is devoted to the study of invertible frame multipliers. These
results are not only of interest in their own right, but also they
pave the way for obtaining some new results for Gabor multipliers
which have been studied mostly by Hans Georg Feichtinger and his
coauthors in recent years.\\


\end{abstract}


\section{\large\bf Introduction}

Throughout this paper, we denote by $\cal H$ a separable Hilbert
space with the inner product ``$\big<\cdot,\cdot\big>"$. Also,
$\ell^2$ and $\ell^\infty$ have their usual meanings and
$(\delta_n)_n$ refers to the canonical orthonormal basis of
$\ell^2$. Moreover, our notation and terminology are standard and,
concerning frames in Hilbert spaces, they are in general those of
the book \cite{c} of Christensen. All over in this paper $\Phi$
and $\Psi$ are sequences $(\varphi_n)_n$ and $(\psi_n)_n$ in $\cal
H$, and in the case where $\Phi$ is a Bessel sequence, the
analysis operator is denoted by $U_\Phi$, the synthesis operator
by $T_\Phi$ and the frame operator by $S_\Phi$. The canonical dual
of the frame $\Phi$ is denoted by
$\widetilde{\Phi}=(\widetilde{\varphi}_n)_n$, and we denote by
$A_{\Phi}$ and $B_{\Phi}$ the lower and upper frame bounds of
$\Phi$, respectively. The notation $m$ is used to denote a complex
scalar sequence $(m_n)_n$, $1/m=(1/m_n)_n$ and
$\overline{m}=(\overline{m}_n)_n$, where ${\overline{m}}_n$
denotes the complex conjugate of $m_n$. The sequence $m$ is called
semi-normalized if $0<\inf_n|m_n|\leq\sup_n|m_n|<\infty$. For
$m\in \ell^\infty$, ${\cal M}_m$ denotes the mapping defined by
${\cal M}_m(c_n)_n=(m_nc_n)_n$ from $\ell^2$ into $\ell^2$. Also,
the main object of study of this work is the operator
$M_{m,\Phi,\Psi}$ which denotes the map defined by the equality
$$M_{m,\Phi,\Psi}(f)=T_\Phi{\cal M}_m U_\Psi(f)=
\sum_{n=1}^\infty m_n\big<f,\psi_n\big>\phi_n
\quad\quad\quad\quad\quad\quad(f\in {\cal H});$$ this operator is
called multiplier with symbol $m$. In particular, the adjoint of
$M_{m,\Phi,\Psi}$, denoted by $(M_{m,\Phi,\Psi})^*$, is equal to
$M_{\overline{m},\Psi,\Phi}$.

The notions of Bessel multiplier, frame multiplier and Riesz
multiplier, as an extension of Gabor multipliers
\cite{feichmulti1}, were introduced and first studied by Balazs
\cite{balaz3} for Hilbert space. As far as we know the subject,
the starting point of the study of such operators is Schatten's
paper \cite{scha}, and in the works
\cite{feichmulti0,feichmulti,feichmulti1}, they are considered in
Fourier and Gabor analysis. In particular, this class of operators
has been extensively studied and has many applications in
different contexts. The reader can find in the papers [1--3] and
[5--12] a lot of information about the history of this class of
operators, some of their properties and their applications in
scientific disciplines and in modern life.

For some applications it is important to consider the stability of
the operator properties of $M_{m,\Phi,\Psi}$ under perturbations
of $\Phi$, $\Psi$ and $m$. For this purpose, we restrict our
investigation to the study of the effect of perturbations of
$\Phi$, $\Psi$ and $m$ on the operator properties of
$M_{m,\Phi,\Psi}$. Among other things, we obtain some conditions
under which the inverse of an invertible frame multiplier can be
represented as a multiplier with the reciprocal symbol and
canonical dual frames of the given ones.


\section{\large\bf Main results}

Let us commence with the following result, which provides some
equivalent conditions for those invertible frame multipliers
$M_{m,\Phi,\Psi}$ whose inverses is
$M_{1/m,\widetilde{\Psi},\widetilde{\Phi}}$.

\begin{theorem}
Suppose that $\Phi$ and $\Psi$ are frames for $\cal H$, and that
$m$ is a semi-normalized sequence for which $M_{m,\Phi,\Psi}$ is
invertible. Suppose also that $A_\Phi$ and $A_\Psi$ are the
optimal lower frame bounds of $\Phi$ and $\Psi$, respectively, and
$|m|$ refers to the sequence $(|m_n|)_n$. Then
$M^{-1}_{m,\Phi,\Psi}=M_{1/m,\widetilde{\Psi},\widetilde{\Phi}}$
if and only if one of the following conditions is satisfied:
\newcounter{j111}
\begin{list}%
{\rm(\roman{j111})}{\usecounter{j111}} \item
$\|S^{-1}_{\Psi}\|=\|M_{m,\Phi,\Psi}^{-1}\,T_\Phi{\cal
M}_{|m|}\|^2$. \item The optimal upper frame bound of the frame
$(M_{m,\Phi,\Psi}^{-1}(m_n\varphi_n))_n$ is $A^{-1}_\Psi$. \item
$\|S^{-1}_{\Phi}\|=\|(M_{m,\Phi,\Psi}^{-1})^*T_\Psi{\cal
M}_{|m|}\|^2$. \item The optimal upper frame bound of the frame
$(M_{\overline{m},\Psi,\Phi}^{-1}(\overline{m}_n\psi_n))_n$ is
$A^{-1}_\Phi$.
\end{list}
\end{theorem}
{\noindent Proof.} Denote $M:=M_{m,\Phi,\Psi}$ and suppose that
$$\psi_n^\dagger=M^{-1}(m_n\varphi_n)\quad\quad
{\hbox{and}}\quad\quad
\varphi_n^\dagger=(M^{-1})^*(\overline{m}_n\psi_n)
\quad\quad\quad\quad\quad\quad(n\in{\Bbb N}).$$ As was shown in
\cite[Theorem 1.1]{balaz3}, $\Psi^\dagger=(\psi_n^\dagger)$ and
$\Phi^\dagger=(\varphi_n^\dagger)$ are the unique dual frames of
$\Psi$ and $\Phi$, respectively, such that
$$M^{-1}=M_{1/m,\Psi^\dagger,\widetilde{\Phi}}\quad\quad\quad{\hbox{and}}
\quad\quad\quad M^{-1}=M_{1/m,\widetilde{\Psi},\Phi^\dagger}.$$
First observe that
\begin{eqnarray*}
S_{\Psi^\dagger}=M^{-1}T_\Phi{\cal M}_{|m|}(M^{-1}T_\Phi{\cal
M}_{|m|})^*~~~~\quad\quad\quad{\hbox{and}}\quad\quad\quad
S_{\widetilde{\Psi}}=S^{-1}_{\Psi},
\end{eqnarray*}
and
\begin{eqnarray*}
S_{\Phi^\dagger}=(M^{-1})^*T_\Psi{\cal
M}_{|m|}{\Big(}(M^{-1})^*T_\Psi{\cal
M}_{|m|}{\Big)}^*\quad\quad{\hbox{and}}\quad\quad\quad
S_{\widetilde{\Phi}}=S^{-1}_{\Phi}.
\end{eqnarray*}
From these, by \cite[Theorem 2.3.1]{mur}, we have
$$\|S_{\Psi^\dagger}\|=\|M^{-1}T_\Phi{\cal M}_{|m|}\|^2\quad\quad
{\hbox{and}}\quad\quad
\|S_{\Phi^\dagger}\|=\|(M^{-1})^*T_\Psi{\cal M}_{|m|}\|^2.$$
Hence, if $\Psi^\dagger=\widetilde{\Psi}$ [respectively,
$\Phi^\dagger=\widetilde{\Phi}$], then the condition (i)
[respectively, (iii)] is satisfied.

Conversely, if condition (i) is satisfied, then, by using
\cite[Theorem 5.7.4]{c} and its proof, there exists a bounded
operator $W:\ell^2\rightarrow{\cal H}$ such that
\begin{equation}\label{parsaj}
\psi_n^\dagger=\widetilde{\psi}_n+V\delta_n\quad\quad\quad\quad\quad\quad(n\in{\Bbb
N}),
\end{equation}
where $V=W(Id_{\cal H}-U_\Psi S^{-1}_\Psi T_{\Psi})$. Now, it is
not hard to check that $T_{\Psi^\dagger}=T_{\widetilde{\Psi}}+V$,
$U_{\Psi^\dagger}=U_{\widetilde{\Psi}}+V^*$,
$T_{\widetilde{\Psi}}V^*=0$ and $VU_{\widetilde{\Psi}}=0$. It
follows that $S_{\Psi^\dagger}=S_{\widetilde{\Psi}}+VV^*$. This,
together with the positivity of the operators $S_{\Psi^\dagger}$,
$S_{\widetilde{\Psi}}$ and $VV^*$ imply that
\begin{eqnarray}\label{gelarj}
\|S_{\Psi^\dagger}\|=\sup_{\|f\|\leq
1}\big<S_{\Psi^\dagger}f,f\big>=\|S_{\widetilde{\Psi}}\|+\sup_{\|f\|\leq
1}{\|}V^*f{\|}^2.
\end{eqnarray}
Therefore, by Eq. (\ref{parsaj}) and (\ref{gelarj}), we get $V=0$
and thus $\Psi^\dagger=\widetilde{\Psi}$. Similarly, if condition
(ii) is satisfied, the inverse of $M$ is
$M_{1/m,\widetilde{\Psi},\widetilde{\Phi}}$ either.

Finally, to prove condition (ii) [respectively, (iv)] is
equivalent to the equality
$M^{-1}=M_{1/m,\widetilde{\Psi},\widetilde{\Phi}}$, it will be
enough to note that, by \cite[Lemma 5.1.6 and Proposition
5.3.8]{c}, $\|S_{\widetilde{\Psi}}\|=A^{-1}_\Psi$ [respectively,
$\|S_{\widetilde{\Phi}}\|=A^{-1}_\Phi$]; this is because of, as
seen above $\|S_{\Psi^\dagger}\|=\|S_{\widetilde{\Psi}}\|$
[respectively, $\|S_{\Phi^\dagger}\|=\|S_{\widetilde{\Phi}}\|$] if
and only if $\Psi^\dagger=\widetilde{\Psi}$ [respectively,
$\Phi^\dagger=\widetilde{\Phi}$], and on the other hand
$\|S_{\Psi^\dagger}\|$ [respectively, $\|S_{\Phi^\dagger}\|$] is
equal to the
optimal upper frame bound of ${\Psi^\dagger}$ [respectively, ${\Phi^\dagger}$].$\e$\\


For the formulation of the following statements, which guarantee
the stability of the operator properties of a frame multiplier
under the perturbations of frames, we need the following
definition. Some basic properties of multipliers can
be found in \cite[Theorem 6.1]{balaz3} and \cite{feichmulti1}.

\begin{definition}
Let $\Phi$ be a sequence in $\cal H$, $m\in\ell^\infty$ and $\mu,
\varepsilon>0$.
\newcounter{j12}
\begin{list}%
{\rm(\roman{j12})}{\usecounter{j12}} \item We say that a sequence
$\Phi'=(\varphi'_n)_n$ in $\cal H$ is a $\mu$-perturbation of
$\Phi$ if $\|T_\Phi-T_{\Phi'}\|\leq\mu$. \item We call the
sequence $m'=(m_n)_n$ a $\varepsilon$-perturbation of $m$ whenever
$\|m-m'\|_\infty\leq\varepsilon$.
\end{list}
\end{definition}

In what follows, for closed subspace $X$ of $\ell^2$, the notation
$\pi_X$ is used to denote the orthogonal projection of $\ell^2$
onto $X$. Moreover, the range of the operator $Q$ is denoted by
${\cal R}(Q)$.

\begin{theorem}\label{per1}
Let $\Phi$ and $\Psi$ be frames for $\cal H$ with frame bounds
$A_\Phi, B_\Phi$ and $A_\Psi, B_\Psi$, respectively, and let $m$
be a semi-normalized symbol. If $\Phi'$ is a $\mu$-perturbation of
$\Phi$ which $\mu<\sqrt{A_\Phi}$, then there exists a frame
$\Psi'=(\psi'_n)_n$ which is a $\lambda\mu$-perturbation of $\Psi$
for some $\lambda>0$ and $M_{m,\Phi',\Psi'}=M_{m,\Phi,\Psi}$. In
particular, the operator properties of $M_{m,\Phi,\Psi}$
{\rm(}such as compactness, invertibility, surjectivity and
etc.{\rm)} are stable under the perturbations of $\Phi$.
\end{theorem}
{\noindent Proof.} First note that \cite[Theorem 5.6.1]{c}
together with the fact that $\Phi'$ is a $\mu$-perturbation of
$\Phi$, where $\mu<\sqrt{A_\Phi}$, implies that $\Phi'$ is a frame
for $\cal H$ with lower frame bound
$A_{\Phi'}=(\sqrt{A_\Phi}-\mu)^2$. Hence, since $m$ is a
semi-normalized sequence, we can deduce that $m\Phi'$ is a frame
for $\cal H$ with lower frame bound
$A_{m\Phi'}=(\inf_n|m_n|)^2(\sqrt{A_\Phi}-\mu)^2$. In particular,
it is easy to see that
$$\ell^2={\cal R}(U_{m\Phi'})\oplus\ker(T_{m\Phi'})\quad\quad\quad
{\hbox{and}}\quad\quad\quad
U_{m\Phi'}T_{\widetilde{m\Phi'}}=U_{m\Phi'}S^{-1}_{m\Phi'}T_{m\Phi'}=\pi_{{\cal
R}(T_{m\Phi'})}.$$ Now, if we set
$$\Psi':={\Big(}M_{\overline{m},\Psi,\Phi}\,S^{-1}_{m\Phi'}(m_n\varphi'_n)+
T_\Psi\pi_{\ker(T_{m\Phi'})}(\delta_n){\Big)}_n\quad\quad\quad(n\in{\Bbb
N}),$$ then, for each sequence $c=(c_n)\in\ell^2$, we have
\begin{eqnarray*}
T_{\Psi'}c&=&\sum_{n=1}^\infty c_nM_{\overline{m},
\Psi,\Phi}\,S^{-1}_{m\Phi'}(m_n\varphi'_n)+\sum_{n=1}^\infty c_n
T_\Psi\pi_{\ker(T_{m\Phi'})}(\delta_n)\\
&=&T_\Psi
U_{m\Phi}T_{\widetilde{m\Phi'}}c+T_\Psi\pi_{\ker(T_{m\Phi'})}c.
\end{eqnarray*}
Therefore, we observe
\begin{eqnarray*}
T_\Psi(U_{m\Phi}-U_{m\Phi'})T_{\widetilde{m\Phi'}}&=&T_\Psi
U_{m\Phi}T_{\widetilde{m\Phi'}}-T_\Psi
U_{m\Phi'}T_{\widetilde{m\Phi'}}+T_\Psi\pi_{\ker(T_{m\Phi'})}-
T_\Psi\pi_{\ker(T_{m\Phi'})}\\
&=&T_\Psi
U_{m\Phi}T_{\widetilde{m\Phi'}}+T_\Psi\pi_{\ker(T_{m\Phi'})}-
(T_\Psi\pi_{{\cal R}(T_{m\Phi'})}+T_\Psi\pi_{\ker(T_{m\Phi'})})\\
&=&T_{\Psi'}-T_{\Psi}.
\end{eqnarray*}
It follows that
\begin{eqnarray*}
\|T_\Psi'-T_\Psi\|&\leq&\|T_\Psi\|\,\|U_{m\Phi}-U_{m\Phi'}\|
\,\|T_{\widetilde{m\Phi'}}\|\\
&\leq&\|T_\Psi\|\,\|m\|_\infty\,\|U_{\Phi}-U_{\Phi'}\|
\,\|T_{\widetilde{m\Phi'}}\|\\
&\leq&\mu\frac{\|m\|_\infty\sqrt{B_\Psi}}{(\inf_n|m_n|)(\sqrt{A_\Phi}-\mu)},
\end{eqnarray*}
and thus $\Psi'$ is a $\lambda\mu$-perturbation of $\Psi$, where
$\lambda:=\|m\|_\infty\sqrt{B_\Psi}/(\inf_n|m_n|)(\sqrt{A_\Phi}-\mu)$.
Finally, we note that
\begin{eqnarray*}
M_{m,\Phi',\Psi'}(f)&=&\sum_{n=1}^\infty\big<f,M_{\overline{m},
\Psi,\Phi}S^{-1}_{m\Phi'}(m_n\varphi'_n)\big>m_n\varphi'_n+
\sum_{n=1}^\infty\big<f,T_\Psi\pi_{\ker(T_{m\Phi'})}(\delta_n)\big>m_n\varphi'_n\\
&=&M_{m,\Phi,\Psi}(f),
\end{eqnarray*}
for all $f\in{\cal H}$. We have now completed the proof of the
theorem.$\e$\\

The following remark is now immediate:

\begin{remark}
{\rm Let $\Phi$ and $\Psi$ be frames for $\cal H$ with frame
bounds $A_\Phi, B_\Phi$ and $A_\Psi, B_\Psi$, respectively, and
let $m$ be a semi-normalized symbol. With an argument similar to
the proof of Theorem \ref{per1} and using the adjoint of the frame
multiplier $M_{m,\Phi,\Psi}$ one can show that if $\Psi'$ is a
$\mu$-perturbation of $\Psi$ which $\mu<\sqrt{A_\Psi}$, then there
exists a frame $\Phi'$ which is a $\lambda\mu$-perturbation of
$\Phi$ for some $\lambda>0$ and
$M_{m,\Phi',\Psi'}=M_{m,\Phi,\Psi}$. In particular, the operator
properties of $M_{m,\Phi,\Psi}$ are stable under the perturbations
of $\Psi$.}
\end{remark}

It is notable that a $\varepsilon$-perturbation of a sequence
$m\in\ell^\infty$ is not necessarily a semi-normalized sequence
even if $m$ is semi-normalized. In the case where $m$ is a
sequence in $\ell^\infty$ for which the inequality $\inf_n|m_n|>
0$ is not necessarily valid, we have the following result. The
reader will remark that a result similar to the following theorem,
which is stated for $\mu$-perturbation of $\Phi$, can be
formulated for $\mu$-perturbation of $\Psi$, and so the details
are omitted here.

\begin{theorem}\label{per2}
Let $\Phi$ and $\Psi$ be frames for $\cal H$ with frame bounds
$A_\Phi, B_\Phi$ and $A_\Psi, B_\Psi$, respectively, and let $m$
be a sequence in $\ell^\infty$ such that the frame multiplier
$M_{m,\Phi,\Psi}$ is invertible. If $\Phi'$ is a
$\mu$-perturbation of $\Phi$ which
$\mu\|m\|_\infty<(\sqrt{B_\Phi}\,\|M_{m,\Phi,\Psi}^{-1}\|)^{-1}$,
then there exists a frame $\Psi'$ which is a
$\lambda\mu$-perturbation of $\Psi$ for some $\lambda>0$ and
$M_{m,\Phi',\Psi'}=M_{m,\Phi,\Psi}$. In particular, the operator
properties of $M_{m,\Phi,\Psi}$ are stable under the perturbations
of $\Phi$.
\end{theorem}
{\noindent Proof.} First note that, it is not hard to check that
$m\Phi$ is a frame with lower frame bound
$A_{m\Phi}:=(B_\Phi\|M_{m,\Phi,\Psi}^{-1}\|^2)^{-1}$. Moreover, we
observe that
$$\|T_{m\Phi}-T_{m\Phi'}\|\leq\|m\|_\infty\|T_{\Phi}-T_{\Phi'}\|\leq\mu\|m\|_\infty.$$
Hence, \cite[Theorem 5.6.1]{c} implies that $m\Phi'$ is a frame for $\cal H$, and thus
if we set
$$\psi'_n:=M_{\overline{m},\Psi,\Phi}\,S^{-1}_{m\Phi'}(m_n\varphi'_n)+
T_\Psi\pi_{\ker(T_{m\Phi'})}(\delta_n)\quad\quad\quad\quad(n\in {\Bbb N}),$$ then with an argument
similar to the proof of Theorem \ref{per1} one can show that
$\Psi'$ is the desired frame.$\e$\\

Next we turn our attention to the perturbations of the symbol of a
frame multiplier whose proof is omitted for conciseness, since it
can be obtained with an argument similar to the proof of Theorems
\ref{per1} and \ref{per2}.

\begin{theorem}\label{per3}
Let $\Phi$ and $\Psi$ be frames for $\cal H$ with frame bounds
$A_\Phi, B_\Phi$ and $A_\Psi, B_\Psi$, respectively, and let $m'$
be a $\varepsilon$-perturbation of $m\in\ell^\infty$. If either
\newcounter{j222}
\begin{list}%
{\usecounter{j222}} \item $\bullet$ the frame multiplier
$M_{m,\Phi,\Psi}$ is invertible and $\varepsilon
B_\Phi<1/\|M_{m,\Phi,\Psi}^{-1}\|$; \item $\bullet$ or $m$ is
semi-normalized and
$\varepsilon\sqrt{B_\Phi}<(\inf_n|m_n|)\sqrt{B_\Phi}$.
\end{list}
\noindent Then in both cases,
$\Psi':={\Big(}M_{\overline{m},\Psi,\Phi}\,S^{-1}_{m'\Phi}(m'_n\varphi_n)+
T_\Psi\pi_{\ker(T_{m'\Phi})}(\delta_n){\Big)}_n$ is a
$\delta\varepsilon$-perturbation of $\Psi$ for some $\delta>0$
and $M_{m',\Phi,\Psi'}=M_{m,\Phi,\Psi}$. In particular, the
operator properties of $M_{m,\Phi,\Psi}$ are stable under the
perturbations of $m$.
\end{theorem}


The next result gives a new representation for the inverse of any
invertible frame multiplier with semi-normalized symbol. In
particular, this proposition shows that the inverse of invertible
frame multiplier $M_{m,\Phi,\Psi}$ has a decomposition into a sum
of $M_{1/m,\widetilde{\Psi},\widetilde{\Phi}}$ and
$\Gamma^*U_{\widetilde{\Phi}}$, where $\Gamma$ is uniquely
determined.

\begin{proposition}\label{main}
Suppose that $\Phi$ and $\Psi$ are frames for $\cal H$, and that
the symbol $m$ is semi-normalized. If $M_{m,\Phi,\Psi}$ is an
invertible multiplier, then there exists a unique bounded operator
$\Gamma:{\cal H}\rightarrow \ell^2$ such that
$$
M_{m,\Phi,\Psi}^{-1}=M_{1/m,\widetilde{\Psi},{\Phi^d}}+
\Gamma^*U_{\Phi^d},
$$
for all dual frames $\Phi^d=(\varphi_n^d)_n$ of $\Phi$.
\end{proposition}
{\noindent Proof.} Define $\Gamma:{\cal H}\rightarrow\ell^2$ by
\begin{equation}\label{1}
\Gamma(f):=U_\Phi(M_{m,\Phi,\Psi}^{-1})^*(f)-{\cal M}_{1/\overline{m}}U_\Psi
S^{-1}_\Psi(f)\quad\quad\quad(f\in {\cal H}).
\end{equation}
Then it is not hard to check that the operator $\Gamma$ is
bounded. In particular,
$M_{m,\Phi,\Psi}^{-1}T_\Phi=\Gamma^*+S^{-1}_\Psi T_\Psi{\cal
M}_{1/m}$. Using any dual frame $\Phi^d$ of $\Phi$ we get
\begin{equation}\label{2}
M_{m,\Phi,\Psi}^{-1}=S^{-1}_\Psi T_\Psi{\cal M}_{1/m}U_{\Phi^d}+
\Gamma^*U_{\Phi^d}.
\end{equation}
It follows that
$M_{m,\Phi,\Psi}^{-1}=M_{1/m,\widetilde{\Psi},{\Phi^d}}+
\Gamma^*U_{\Phi^d}$ for all dual frames $\Phi^d$ of $\Phi$.
Finally, the proof will be completed by showing that the operator
$\Gamma$ is uniquely determined. To this end, suppose on the
contrary that Eq. (\ref{2}) are hold for two operators $\Gamma_1$
and $\Gamma_2$. Hence, we have $\Gamma_1^*U_{\Phi^d}=
\Gamma_2^*U_{\Phi^d}$ for all dual frames $\Phi^d$ of $\Phi$. We
now invoke part (i) of
\cite[Theorem 1.2]{balaz11} to conclude that $\Gamma_1=\Gamma_2$.$\e$\\

The following remark is now immediate:

\begin{remark}
{\rm Suppose that $\Phi$ and $\Psi$ are frames for $\cal H$, and
that the symbol $m$ is semi-normalized. If $M_{m,\Phi,\Psi}$ is an
invertible multiplier, then
\newcounter{j1}
\begin{list}%
{\rm(\roman{j1})}{\usecounter{j1}} \item For operator $\Gamma$ in
Proposition \ref{main} it is not hard to check that $T_\Psi\Gamma=0$.
It follows that in the case where $\Psi$ is a Riesz basis, then
$M_{m,\Phi,\Psi}^{-1}=M_{1/m,\widetilde{\Psi},{\Phi^d}}$ for all
dual frames $\Phi^d$ of $\Phi$.
 \item  It can be shown by routine
calculations that if $\Psi$ is equivalent to $m\Phi$, then for
each dual frames $\Phi^d$ of $\Phi$ the inverse of
$M_{m,\Phi,\Psi}$ is $M_{1/m,\widetilde{\Psi},\Phi^d}$.
Conversely, if
$M_{m,\Phi,\Psi}^{-1}=M_{1/m,\widetilde{\Psi},\Phi^d}$ for all
dual frames $\Phi^d$ of $\Phi$, then Proposition \ref{main}
implies that $\Gamma^*U_{\Phi^d}=0$ for all dual frames $\Phi^d$
of $\Phi$. This together with \cite[Theorem 1.2(i)]{balaz11}
implies that $\Gamma=0$. From this, by Eq. (\ref{1}), we deduce
that ${\cal M}_{\overline{m}}U_\Phi
(M_{m,\Phi,\Psi}^{-1})^*=U_\Psi S^{-1}_\Psi$. It follows that
$$\big<f,m_n\phi_n\big>=\big<f,M_{m,\Phi,\Psi}S^{-1}_\Psi \psi_n\big>\quad\quad\quad
\quad\quad\quad{\Big(}f\in {\cal H}~ {\hbox{and}}~ n\in {\Bbb
N}{\Big)},$$ and thus the frames $\Psi$ and $m\Phi$ are
equivalent. \item With an argument similar to the proof of
Proposition \ref{main} above one can show that $\Theta=U_\Psi
M^{-1}_{m,\Phi,\Psi}-{\cal M}_{1/m}U_\Phi S^{-1}_\Phi$ is the
unique bounded operator such that
$$
M^{-1}_{m,\Phi,\Psi}=M_{1/m,{\Psi^d},\widetilde{\Phi}}+T_{\Psi^d}\Theta,
$$
for all dual frames $\Psi^d$ of $\Psi$. In particular,
$T_\Phi\Theta=0$. Hence, if $\Phi$ is a Riesz basis, then
$M^{-1}_{m,\Phi,\Psi}=M_{1/m,{\Psi^d},\widetilde{\Phi}}$ for all
dual frames $\Psi^d$ of $\Psi$. Moreover, for each dual frames
$\Psi^d$ of $\Psi$
$M_{m,\Phi,\Psi}^{-1}=M_{1/m,{\Psi^d},\widetilde{\Phi}}$  if and
only if $\Phi$ is equivalent to $\overline{m}\Psi$.\\
\end{list}
}
\end{remark}

\vspace{3mm}

\noindent {\sc Hossein Javanshiri}\\
Department of Mathematics,
Yazd University,
P.O. Box: 89195-741, Yazd, Iran\\
E-mail: h.javanshiri@yazd.ac.ir\\


\footnotesize



\begin{thebibliography}{HD}




\bibitem{balaz3} {\sc P. Balazs},  Basic definition and properties of Bessel
multipliers, {\it J. Math. Anal. Appl.} {\bf 325} (2007), 571-–85.



\bibitem{balaz11} {\sc P. Balazs and D. T. Stoeva}, Representation of
the inverse of a frame multiplier, {\it J. Math. Anal. Appl.}
{\bf 422} (2015), 981--994.


\bibitem{bened} {\sc J. Benedetto and G. Pfander}, Frame expansions for Gabor multipliers,
{\it Appl. Comput. Harmon. Anal.} {\bf 20} (2006), 26--40.


\bibitem{c} {\sc O. Christensen,} Introduction to frames and
Riesz bases, Birkh\"{a}user, (2003).


\bibitem{cord} {\sc E. Cordero and F. Nicola}, Remarks on Fourier multipliers and
applications to the wave equation,
{\it J. Math. Anal. Appl.} {\bf 353} (2009), 583--591.


\bibitem{cord1} {\sc E. Cordero and K. Gr\"{o}chenig},
Necessary conditions for Schatten class localization
operators, {\it Proc. Amer. Math. Soc.} {\bf 133} (2005), 3573--3579


\bibitem{cord2} {\sc E. Cordero, K. Gr\"{o}chenig and F. Nicola},
Approximation of Fourier integral operators by Gabor multipliers,
{\it J. Fourier Anal. Appl.} {\bf 18} (2012), 661--684.



\bibitem{dorf} {\sc
M. D\"{o}rfler and B. Torr\'{e}sani},
Representation of operators in the time-frequency
domain and generalized Gabor multipliers,
{\it J. Fourier Anal. Appl.} {\bf 16} (2010), 261--293.


\bibitem{feichmulti0} {\sc
H. G. Feichtinger, M. Hampejs and G. Kracher,}
Approximation of matrices by Gabor multipliers, {\it IEEE
Signal Process. Lett.} {\bf 11} (2004), 883--886.



\bibitem{feichmulti} {\sc H.G. Feichtinger and G. Narimani},
Fourier multipliers of classical modulation spaces,
{\it Appl. Comput. Harmon. Anal.} {\bf 21} (2006), 349--359.



\bibitem{feichmulti1} {\sc H.G. Feichtinger and K. Nowak},
A first survey of Gabor multipliers. Advances in Gabor analysis.
Appl. Numer. Harmon. Anal., pp. 99--128. Birkh\"{a}user, Boston
(2003).



\bibitem{gro1} {\sc
K. Gr\"{o}chenig}, {Representation and approximation of pseudodifferential
operators by sums of Gabor multipliers,} {\it Appl. Anal.} {\bf 90} (2011), 385--401.



\bibitem{mur} {\sc G. J. Murphy}, C$^*$-Algebras and Operator Theory,
Academic Press, London (1990).



\bibitem{scha} {\sc R. Schatten}, Norm Ideals of Completely
Continuous Operators, Springer, Berlin, (1960).




\end{thebibliography}
\end{document}